\documentstyle[12pt]{amsart}
\newtheorem{theorem}{Theorem}[section]

\newtheorem{remark}[theorem]{Remark}
\newtheorem{lemma}[theorem]{Lemma}
\newtheorem{proposition}[theorem]{Proposition}

\numberwithin{equation}{section}

\begin{document}

\title [Symplectic and Lagrangian mean curvature flows]
 {Translating solitons to symplectic and Lagrangian mean curvature flows}
\author{Xiaoli Han, Jiayu Li}

\address{Math. Group, The abdus salam ICTP\\ Trieste 34100,
   Italy}
\email{xhan@@ictp.it}

\address{Math. Group, The abdus salam ICTP\\ Trieste 34100,
   Italy\\
   and Academy of Mathematics and Systems Sciences\\ Chinese Academy of
Sciences\\ Beijing 100080, P. R. of China. } \email{jyli@@ictp.it}

\keywords{Symplectic surface, Lagrangian surface, translating
slotion, mean curvature flow.}

\date{}

\begin{abstract}
In this paper, we construct finite blow-up examples for symplectic
mean curvature flows and we study properties of symplectic
translating solitons. We prove that, the K\"ahler angle $\alpha$
of a symplectic translating soliton with $\max |A|=1$ satisfies
that $\sup |\alpha|>\frac{\pi}{4}\frac{|T|}{|T|+1}$ where $T$ is
the direction in which the surface translates.
\end{abstract}

\maketitle

 Mathematics Subject Classification (2000): 53C44 (primary),
53C21 (secondary).

\section{Introduction}
We consider the evolution of a symplectic surface (Lagrangian
surface) in a K\"ahler-Einstein surface by its mean curvature,
which we call a symplectic mean curvature flow (Lagrangian mean
curvature flow). We \cite{HL1} showed that, if the scalar
curvature of the K\"ahler-Einstein surface is positive and the
initial surface is sufficiently close to a holomorphic curve, then
the symplectic mean curvature flow exists globally and converges
to a holomorphic curve at infinity. In this paper, we construct
examples to show that, in general the symplectic mean curvature
flow may blow-up at a finite time. Therefore, it is necessary to
study the singularity of a symplectic mean curvature flow.

Chen-Li (\cite{CL1}, \cite{CL2}, \cite{CL3}) and Wang \cite {Wa}
independently proved that there is no Type I singularity along the
symplectic (almost calibrated Lagrangian) mean curvature flow.
Therefore, the structure of the type II singularity has been of
great interest.

One of the most important examples of Type II singularity is the
translating soliton (c.f. \cite {Ha1}, \cite{HS}), which is one
class of eternal solutions defined for $-\infty<t<\infty$. More
precisely, the translating solitons are surfaces which evolve by
translating in space with a constant velocity. Hamilton studied
this kind of eternal solutions to the mean curvature flow of a
hypersurface in ${\mathbf{R}}^n$ \cite{Ha1} and to the Ricci flow
\cite{Ha2}, \cite{Ha3}. The main purpose of this paper is to study
the properties of symplectic translating solitons.

\vspace{.2in}

\noindent {\bf Definition } {\it The translating soliton is called
a standard translating soliton if the norm of its second
fundamental form $A$ satisfies that $\max |A|=1$. Let $$
ST=\{\Sigma~|~\Sigma~{\rm is~a~standard~symplectic~ translating
~soliton }.\}$$}

\vspace{.2in}

\noindent {\bf Remark } {\it Since we are interested in the
translating solitons which arise in the blow-up analysis of
singularities, it is reasonable to assume that $\max |A|=1$. }

\vspace{.2in}

\noindent {\bf Main Theorem 1 } {\it Suppose that $\Sigma\in ST$
is a symplectic standard translating soliton with K\"ahler angle
$\alpha$, then $\sup_{\Sigma} |\alpha|>
\frac{\pi}{4}\frac{|T|}{|T|+1}$, where $T$ is the direction in
which the surface translates.}

\vspace{.2in}

Note that, the "grim reaper" $(x, y, -\ln\cos x, 0),~ |x|<\pi/2,~
y\in \mathbf{R}$ is one standard symplectic translating soliton
which translates in the direction of the constant vector $(0, 0,
1, 0)$. In this case the K\"ahler angle $\alpha= x$ and it is
clear that $\sup_\Sigma |\alpha|=\pi/2$, i.e,
$\inf_{\Sigma}\cos\alpha=0$. In the study of symplectic mean
curvature flows of compact surfaces, we assume that
$\cos\alpha\geq\delta>0$ on the initial surface and consequently
\cite{CL1} at each time $t$, on $\Sigma_t$,
$\cos\alpha\geq\delta>0$. One of the purpose of studying
symplectic translating solitons is to rule out them as the
limiting flows in the blow-up analysis at a type II singularity.
More precisely, we believe,

\vspace{.2in}

\noindent {\bf Conjecture 1 } {A symplectic translating soliton
can not be a limiting flow of rescaled surfaces at a type II
singular point.}

\vspace{.2in}

Even more we conjecture that

\vspace{.2in}

\noindent {\bf Conjecture 2 } {\it Any blowing up at a type II
singularity of a symplectic mean curvature flow is a union of
non-flat minimal surfaces in ${\mathbb{R}}^4$.}

\vspace{.2in}

As an analogous result of translating solitons in almost
calibrated mean curvature flow \cite{Sm1}, we have,

\vspace{.2in}

\noindent {\bf Main Theorem 2 } {\it Suppose that $\Sigma$ is a
standard translating soliton to the almost calibrated Lagrangian
mean curvature flow with Lagrangian angle $\theta$, then
$\sup_\Sigma |\theta|>\frac{\sqrt{2}\pi}{4}\frac{|T|}{|T|+1}$,
where $T$ is the direction in which the surface translates.}
\vspace{.2in}

The authors thank the referee for his many helpful comments.

\vspace{.2in}

\section{Preliminaries}

In this section, we fix some notations and recall some basic facts
on symplectic mean curvature flows and Lagrangian mean curvature
flows. We consider immersions
$$ F_0: \Sigma\rightarrow {\mathbf{R}}^4
$$ of smooth surface $\Sigma$ in ${\mathbf{R}}^4$. If $\Sigma$ evolves
along the mean curvature, then there is a one-parameter family
$F_t=F(\cdot, t)$ of immersions which satisfy the mean curvature
flow equation:
\begin{equation}\left\{
\begin{array}{clcr}\frac{d}{dt}F(x, t)&=&H(x, t)\\ F(x,
0)&=&F_0(x).\end{array}\right.
\end{equation} Here $H(x, t)$ is the mean curvature vector of
$\Sigma_t=F_t(\Sigma)$ at $F(x, t)$.

Let $\omega$ and $\langle\cdot, \cdot\rangle$ denote the standard
K\"ahler form and Euclidean metric on ${\mathbf{R}}^4$
respectively. We choose a local field of orthonormal frames $e_1,
e_2$, $v_1, v_2$ of ${\mathbf{R}}^4$ along $\Sigma$ such that
$e_1, e_2$ are tangent vectors of $\Sigma$ and $v_1, v_2$ are in
the normal bundle over $\Sigma$. Denote the induced metric on
$\Sigma$ by $(g_{ij})$. The second fundamental form $A$ and the
mean curvature vector $H$ of $\Sigma$ can be expressed, in the
local frame, as $ {A}=A^\alpha v_\alpha$, and ${H}=-H^\alpha
v_\alpha$, where and throughout this paper all repeated indices
are summed over suitable range. For each $\alpha$, the coefficient
$A^\alpha$ is a $2\times 2$ matrix $(h^\alpha_{ij})_{2\times 2}$.
By the Weingarten equation (cf. \cite{Sp}), we have
\begin{eqnarray*} h^\alpha_{ij}&=&\langle v_\alpha,
\overline\nabla_i\overline\nabla_j F\rangle=-\langle
\overline\nabla_j v_\alpha, \overline\nabla_i
F\rangle=h^\alpha_{ji},\\
H^\alpha &=&g^{ij}h^\alpha_{ij}=h^\alpha_{ii},
\end{eqnarray*} where $\overline\nabla$ is the connection on
${\mathbf{R}}^4$. The norm of the second fundamental form of
$\Sigma$ is:
$$|{A}|^2=\sum_{\alpha}|A^\alpha|^2=g^{ij}g^{kl}h^\alpha_{ik}
h^\alpha_{jl}=h^\alpha_{ik}h^\alpha_{ik}.
$$ A surface $\Sigma$ is
called {\it symplectic} if the K\"ahler angle $\alpha$ (c.f.
\cite{CW}) satisfies that $ \cos\alpha> 0, $ where $\cos\alpha$ is
defined by $\omega |_\Sigma=\cos\alpha d\mu_\Sigma$, where
$d\mu_{\Sigma}$ is the induced volume form of $\Sigma$, and
$\omega$ is the standard K\"ahler form on ${\mathbf{C}}^2$.

Suppose $\Sigma$ is symplectic and evolves along the mean
curvature in ${\mathbf{R}}^4$ which is called symplectic mean
curvature flow. Let $J_{\Sigma_t}$ be an almost complex structure
in a tubular neighborhood of $\Sigma$ in ${\mathbf{R}}^4$ with
\begin{equation}\label{eq2}
\left\{\begin{array}{clcr} J_{\Sigma_t}e_1&=&e_2\\
J_{\Sigma_t}e_2&=&-e_1\\ J_{\Sigma_t}v_1&=&v_2\\
J_{\Sigma_t}v_2&=&-v_1.
\end{array}\right.
\end{equation}

It is not difficult to verify (\cite{CT} and \cite{CL1}) that,
\begin{eqnarray*}
|\overline\nabla J_{\Sigma_t}|^2&=&
|h_{11}^2+h_{12}^1|^2+|h^2_{21}+h^1_{22}|^2
+|h_{12}^2-h_{11}^1|^2+|h^2_{22}-h^1_{21}|^2 \\
&=&\frac{1}{2}|H|^2+\frac{1}{2}\left(((h^1_{11}+h^1_{22})-2(h^2_{12}+h^1_{22}))^2
+((h^2_{11}+h^2_{22})-2(h^1_{21}+h^2_{11}))^2\right)\\
&\geq& \frac{1}{2}|H|^2.
\end{eqnarray*}

Recall that (\cite{CL1}) the K\"ahler angle $\alpha$ of $\Sigma_t$
in ${\mathbf{R}}^4$ satisfies the parabolic equation:
\begin{equation}\label{alpha} \left(\frac{\partial}{\partial t}-\Delta
\right)\cos\alpha = |\overline\nabla J_{\Sigma_t}|^2\cos\alpha.
\end{equation}

Suppose that the initial surface is symplectic, i.e., $\cos
\alpha> 0$ , then by applying the parabolic maximum principle to
this evolution equation, one concludes that $\cos\alpha$ remains
positive as long as the mean curvature flow has a smooth solution
(cf. \cite{CT1}, \cite{CL1}, \cite{Wa}).

Let $J$ denote the standard complex structure on ${\mathbf{C}}^2$.
We also consider a parallel holomorphic $(2, 0)$ form,
$$ \Omega=dz_1\wedge dz_2. $$
 A surface $\Sigma$ is said to be {\it Lagrangian} if
$\omega |_\Sigma=0$. This implies that (see \cite{HaL})
$$\Omega |_\Sigma=e^{i\theta}d\mu_\Sigma,
$$
where $d\mu_\Sigma$ denotes the induced volume form of $\Sigma$
and $\theta$ is called Lagrangian angle which is some multivalued
function. If $ \cos\theta> 0, $ then $\Sigma$ is called {\it
almost-calibrated}. The relation between the Lagrangian angle and
the mean curvature vector is given in \cite{HaL} (also see
\cite{TY}),
$$H=J\nabla\theta.$$

Suppose $\Sigma$ is Lagrangian and evolves by the mean curvature,
Smoczyk has shown that (\cite{Sm1}, \cite{Sm2}, \cite{Sm3}),
\begin{equation}\label{theta} (\frac{\partial}{\partial
t}-\Delta)\cos\theta=|H|^2\cos\theta.
\end{equation} If the initial surface is almost calibrated, then for every
$t$, $\Sigma_t$ is also almost calibrated, i.e. $\cos\theta> 0$,
along the mean curvature flow by the parabolic maximum principle.

\section{Property of translating solitons}

Suppose that $\Sigma_t$ is a translating soliton which translates
in the direction of the constant vector $T$. That means
$F_t=F+tT$, i.e, $\Sigma_t=\Sigma+tT$. Let $V=v^i e_i$ be the
tangent part of $T$. Then the normal component must be $N=H^\alpha
v_\alpha$ to solve the mean curvature flow. If we take the
equation $v^i e_i+H^\alpha v_\alpha=T$ and differentiate it, then
we get
$$\overline\nabla_j v^i e_i+v^i h^\alpha_{ij}v_\alpha+\overline\nabla_j
H^\alpha v_\alpha-H^\alpha h^\alpha_{ij} e_i=0,
$$ where we assume that $\nabla_{e_i}e_j=0$ at the considered point.
Separating the tangential and normal part we get that,
\begin{eqnarray}\label{T}
\overline\nabla_j v^i&=& H^\alpha h^\alpha_{ij}\nonumber\\
\overline\nabla_j H^\alpha &=&-v^i h^\alpha_{ij}.
\end{eqnarray}
Using these equations we can get the following identities for the
translating solitons.

\begin{proposition}\label{p1} On the translating soliton, the K\"ahler
angle and the Lagrangian angle satisfy the following elliptic
equations,
\begin{equation}\label{e1} -\Delta\cos\alpha=|\overline\nabla J_{\Sigma_t}|^2\cos\alpha+v^i
\nabla_i\cos\alpha, \end{equation} and
\begin{equation}\label{e2}
 -\Delta\cos\theta=|H|^2\cos\theta+v^i
\nabla_i\cos\theta
\end{equation}
\end{proposition}

{\it Proof.} Along the translating soliton, we have
$T(\cos\alpha)=0$. It implies that $H^\alpha
v_\alpha(\cos\alpha)=-v^i \nabla_i\cos\alpha$. In (\ref{alpha}),
we have  $\frac{\partial}{\partial t}\cos\alpha=H^\alpha
v_\alpha(\cos\alpha)$, which implies that
\begin{eqnarray*}
-\Delta\cos\alpha &=&|\overline\nabla
J_{\Sigma_t}|^2\cos\alpha+v^i \nabla_i\cos\alpha.
\end{eqnarray*} Using (\ref{theta}), similarly we can show (\ref{e2}).
 \hfill Q. E. D.

By solving the system of equations,
\begin{eqnarray*} \overline\nabla_j H^\alpha &=&-v^i h^\alpha_{ij},\\
h^1_{11}+h^1_{22}&=&H^1,\\ h^2_{11}+h^2_{22}&=&H^2,\end{eqnarray*}
we get that at the points $|V|\neq 0$,
\begin{eqnarray*}
h^1_{11} &=& -(v^1 \overline\nabla_1 H^1-v^2\overline\nabla_2
H^1-H^1(v^2)^2)/|V|^2,\\h^1_{12} &=& -(v^1 \overline\nabla_2
H^1+v^2\overline\nabla_1 H^1+H^1v^2v^1)/|V|^2,\\h^1_{22} &=& (v^1
\overline\nabla_1 H^1-v^2\overline\nabla_2
H^1+H^1(v^1)^2)/|V|^2,\\h^2_{11} &=& -(v^1 \overline\nabla_1
H^2-v^2\overline\nabla_2 H^2-H^2(v^2)^2)/|V|^2,\\h^2_{12} &=&
-(v^1 \overline\nabla_2 H^2+v^2\overline\nabla_1
H^2+H^2v^2v^1)/|V|^2,\\h^2_{22} &=& (v^1 \overline\nabla_1
H^2-v^2\overline\nabla_2 H^2+H^2(v^1)^2)/|V|^2.
\end{eqnarray*} Then it is easy to check that,

\begin{proposition}\label{A} On the translating soliton, at the
points $|V|\neq 0$,
$$ |A|^2=|H|^2+2\frac{|\nabla H|^2}{|V|^2}+\frac{v^i\nabla_i
|H|^2}{|V|^2}. $$
\end{proposition}

\begin{proposition} \label{p2} On the translating soliton, at the points $|V|\neq 0$,
the mean curvature vector satisfies that,
\begin{eqnarray}\label{e3}\Delta |H|^2 &=& 2|\nabla H|^2-2(H^\alpha
h^\alpha_{ij})^2-v^i\nabla_i |H|^2\nonumber\\ &=& 2|\nabla
H|^2-2|H|^4-v^i \nabla_i |H|^2-\frac{|\nabla
|H|^2|^2}{|V|^2}-2\frac{|H|^2}{|V|^2}v^i \nabla_i
|H|^2.\end{eqnarray}
\end{proposition}

{\it Proof.} Using equation (\ref{T}) again, we have,
\begin{eqnarray*} \Delta |H|^2&=&2H^\alpha\Delta
H^\alpha+2|\nabla H^\alpha|^2\\ &=&-2H^\alpha \overline\nabla_j(v^i h^\alpha_{ij})+2|\nabla H|^2 \\
&=&-2H^\alpha v^i\overline\nabla_i H^\alpha-2H^\alpha
h^\alpha_{ij}\overline\nabla_j v^i+2|\nabla H|^2\\ &=& 2|\nabla
H|^2-v^i\nabla_i |H|^2-2(H^\alpha h^\alpha_{ij})^2.
\end{eqnarray*}
Putting $h^\alpha_{ij}$ into the above equation, we can prove the
proposition. \hfill Q. E. D.

In the sequence we prove one monotonicity formula for translating
solitons. Let $H(X, X_0, t_0, t)$ be the backward heat kernel on
${\mathbf{R}}^4$. Define
\begin{eqnarray*}
\rho (F, t)&=& 4\pi(t_0-t)H(F, X_0, t_0, t)\\ &=&
\frac{1}{4\pi(t_0-t)}\exp(-\frac{|F-X_0|^2}{4(t_0-t)})
\end{eqnarray*} for $t<t_0$.

\begin{proposition}\label{p3}
On the translating soliton, the K\"ahler angle satisfies the
following formula,
\begin{eqnarray}\label{mono1}
&&\frac{\partial}{\partial
t}\left(\int_{\Sigma_t}\frac{1}{\cos\alpha} \rho(F, X_0, t,
t_0)d\mu_t\right)\nonumber\\&&= -\left(
\int_{\Sigma_t}\frac{1}{\cos\alpha}\rho (F, X_0, t, t_0)
\left|H+\frac{(F-X_0)^{\perp}}{2(t_0-t)}\right|^2d\mu_t \right.\nonumber \\
&&\left.+ \int_{\Sigma_t}\frac{1}{\cos\alpha}\rho (F, X_0, t,
t_0)\left|\overline{\nabla}J_{\Sigma_t}\right|^2d\mu_t
\right.\nonumber\\
&&\left.+\int_{\Sigma_t}\frac{2}{\cos^3\alpha}\left|\nabla
\cos\alpha\right|^2\rho (F, X_0, t, t_0)d\mu_t\right).
\end{eqnarray}
\end{proposition}

{\it Proof.}  A straightforward calculation shows that,
\begin{eqnarray*}
\frac{\partial}{\partial t}\rho &=&
(\frac{1}{t_0-t}-\frac{\frac{\partial F}{\partial t}\cdot
(F-X_0)}{2(t_0-t)}-\frac{|F-X_0|^2}{4(t_0-t)^2} )\rho\\ &=&
(\frac{1}{t_0-t}-\frac{T\cdot
(F-X_0)}{2(t_0-t)}-\frac{|F-X_0|^2}{4(t_0-t)^2} )\rho \\
&=& (\frac{1}{t_0-t}-\frac{H\cdot (F-X_0)}{2(t_0-t)}-\frac{V\cdot
(F-X_0)}{2(t_0-t)}-\frac{|F-X_0|^2}{4(t_0-t)^2})\rho.
\end{eqnarray*} On $\Sigma_t$ we have,
\begin{eqnarray*}
\Delta\rho=(\frac{\langle F-X_0, \nabla
F\rangle^2}{4(t_0-t)^2}-\frac{\langle F-X_0, \Delta
F\rangle}{2(t_0-t)}-\frac{|\nabla F|^2}{2(t_0-t)})\rho,
\end{eqnarray*}
where $\nabla, \Delta$ are connection and Laplacian on $\Sigma_t$
in the induced metric. Note that,
$$ |\nabla F|^2=2,~~~~~~~~~~~~~\Delta F=H,
$$ and add these two equations together, we get,
\begin{eqnarray*}
(\frac{\partial}{\partial t}+\Delta)\rho=-\left(\left
|H+\frac{(F-X_0)^\perp}{2(t_0-t)}\right|^2-|H|^2+\frac{V\cdot
(F-X_0)}{2(t_0-t)}\right)\rho.
\end{eqnarray*}
Now we have, \allowdisplaybreaks \begin{eqnarray*}&&
\frac{\partial}{\partial t}\int_{\Sigma_t}\frac{1}{\cos\alpha}
\rho(F, X_0, t, t_0)d\mu_t
\\ &&=\int_{\Sigma_t}\frac{1}{\cos\alpha}(\frac{\partial}{\partial
t}+\Delta)\rho d\mu_t-\int_{\Sigma}\frac{1}{\cos\alpha}\Delta\rho
d\mu_t\\
&&=
-\int_{\Sigma_t}\frac{1}{\cos\alpha}\left(|H+\frac{(F-X_0)^\perp}{2(t_0-t)}|^2
-|H|^2+\frac{V\cdot (F-X_0)}{2(t_0-t)}\right)\rho
d\mu_t\\&&-\int_{\Sigma_t}\Delta(\frac{1}{\cos\alpha})\rho d\mu_t.
\end{eqnarray*} Using equation (\ref{e1}), we get that,
\begin{eqnarray}\label{e4}
&& \frac{\partial}{\partial t}\int_{\Sigma_t}\frac{1}{\cos\alpha}
\rho(F, X_0, t, t_0)d\mu_t\nonumber\\
&&=-\int_{\Sigma_t}\frac{1}{\cos\alpha}\left(|H+\frac{(F-X_0)^\perp}{2(t_0-t)}|^2
-|H|^2\right)\rho d\mu_t+\int_{\Sigma_t}\frac{1}{\cos\alpha}
V(\rho) d\mu_t\nonumber\\&&-\int_{\Sigma_t}
(\frac{|\overline\nabla J_{\Sigma_t}|^2}{\cos\alpha}-V
(\frac{1}{\cos\alpha}))\rho
d\mu_t-\int_{\Sigma_t}\frac{2}{\cos^3\alpha}|\nabla\cos\alpha|^2\rho
d\mu_t\nonumber\\
&&=-\int_{\Sigma_t}\frac{1}{\cos\alpha}\left(|H+\frac{(F-X_0)^\perp}{2(t_0-t)}|^2
-|H|^2\right)\rho d\mu_t\nonumber\\
&&-\int_{\Sigma_t}\frac{2}{\cos^3\alpha}|\nabla\cos\alpha|^2\rho
d\mu_t-\int_{\Sigma_t} \frac{|\overline\nabla
J_{\Sigma_t}|^2}{\cos\alpha}\rho d\mu_t\nonumber
\\&&+\int_{\Sigma_t}  V (\frac{1}{\cos\alpha}\rho
d\mu_t)-\int_{\Sigma_t} \frac{1}{\cos\alpha}\rho
V(d\mu_t)\nonumber\\&&=-\int_{\Sigma_t}\frac{1}{\cos\alpha}\left(|H+\frac{(F-X_0)^\perp}{2(t_0-t)}|^2
-|H|^2\right)\rho d\mu_t
\nonumber\\&&-\int_{\Sigma_t}\frac{2}{\cos^3\alpha}|\nabla\cos\alpha|^2\rho
d\mu_t-\int_{\Sigma_t} \frac{|\overline\nabla
J_{\Sigma_t}|^2}{\cos\alpha}\rho d\mu_t \nonumber
\\&&-\int_{\Sigma_t} H(\frac{1}{\cos\alpha}\rho
d\mu_t)+\int_{\Sigma_t} \frac{1}{\cos\alpha}\rho H(d\mu_t).
\end{eqnarray} Recall that
$$ H(d\mu_t)=\frac{\partial}{\partial t}(d\mu_t)=-|H|^2 d\mu_t,
$$ thus,
$$\int_{\Sigma_t} \frac{1}{\cos\alpha}\rho H(d\mu_t)
=-\int_{\Sigma_t} \frac{|H|^2}{\cos\alpha}\rho d\mu_t.
$$ Since $H$ is normal to $\Sigma_t$, we have,
\begin{eqnarray*} \int_{\Sigma_t} H(\frac{1}{\cos\alpha}\rho d\mu_t) &=&
H(\int_{\Sigma_t} \frac{1}{\cos\alpha}\rho d\mu_t)\\
&=&-V(\int_{\Sigma_t} \frac{1}{\cos\alpha}\rho d\mu_t)=0.
\end{eqnarray*}

Adding these identities into (\ref{e4}), we can prove the
proposition. \hfill Q. E. D.

By the same argument, we can also get one monotonicity formula for
Lagrangian angle.
\begin{proposition}\label{p4}
On the translating soliton, the Lagrangian angle satisfies the
following formula,
\begin{eqnarray}\label{mono2}
&&\frac{\partial}{\partial
t}\left(\int_{\Sigma_t}\frac{1}{\cos\theta} \rho(F, X_0, t,
t_0)d\mu_t\right)\nonumber\\&&= -\left(
\int_{\Sigma_t}\frac{1}{\cos\theta}\rho (F, X_0, t, t_0)
\left|H+\frac{(F-X_0)^{\perp}}{2(t_0-t)}\right|^2d\mu_t \right.\nonumber \\
&&\left.+ \int_{\Sigma_t}\frac{1}{\cos\theta}\rho (F, X_0, t,
t_0)\left|H\right|^2d\mu_t\right.\nonumber\\
&&\left.+\int_{\Sigma_t}\frac{2}{\cos^3\theta}\left|\nabla
\cos\theta\right|^2\phi\rho (F, X_0, t, t_0)d\mu_t\right).
\end{eqnarray}
\end{proposition}

\section {Proof of Main theorems}

Using the gradient estimates, we can prove our main theorem. First
we give an estimate of the evolution equation of $|H|^2$. Recall
that we assume $|A|^2\leq 1$. Applying the equations in
Proposition \ref{p2}, we have that if $ |H|^2\geq\varepsilon$,
then
\begin{eqnarray*}
\Delta |H|^2 &\geq& 2|\nabla H|^2-2|H|^2-v^i\nabla_i |H|^2\\
&\geq& 2|\nabla H|^2-\frac{2}{\varepsilon}|H|^4-v^i \nabla_i
|H|^2,
\end{eqnarray*} and if $|H|^2<\varepsilon$, then by Cauchy-Schwartz
inequality we have,
\begin{eqnarray*}
\Delta |H|^2 &=& 2|\nabla H|^2-2|H|^4-v^i \nabla_i
|H|^2-\frac{|\nabla |H|^2|^2}{|V|^2}-2\frac{|H|^2}{|V|^2}v^i
\nabla_i |H|^2 \\ &\geq& 2|\nabla H|^2-(2+\delta)|H|^4-v^i
\nabla_i |H|^2-4(1+\frac{1}{\delta})\frac{\varepsilon
}{|T|-\varepsilon} |\nabla |H||^2.
\end{eqnarray*} Set $\tilde a=\frac{2}{\varepsilon}=2+\delta$,
$\tilde b=\frac{1}{2}+ (1+\frac{1}{\delta})\frac{\varepsilon
}{|T|-\varepsilon},$ then in any case we have,
\begin{eqnarray}\label{e8}
\Delta |H|^2 \geq (4-4\tilde b)|\nabla |H||^2-\tilde
a|H|^4-v^i\nabla_i |H|^2.
\end{eqnarray}

We first prove Main Theorem in Lagrangian case.

\vspace{.1in}

\noindent {\bf Main Theorem 2 } {\it Suppose that $\Sigma$ is a
standard translating soliton to the almost calibrated Lagrangian
mean curvature flow with Lagrangian angle $\theta$, then
$\sup_\Sigma |\theta|>\frac{\sqrt{2}\pi}{4}\frac{|T|}{|T|+1}$,
where $T$ is the direction in which the surface translates.}

\vspace{.1in}

{\it Proof.} We argue it by contradiction. Suppose there is a
standard translating soliton with
$\cos(\sqrt{ab}\theta)\geq\delta_0$ for some $\delta_0>0$, where
$a, b$ are constants which are determined later. Now we consider
the function
$$f=\frac{|H|^2}{\cos^{\frac{1}{b}}(\sqrt{ab}\theta)}.$$ It is well
know that
$$ -\Delta\theta=H^\alpha v^\alpha(\theta)=-v^i\nabla_i \theta.$$
Using (\ref{e2}) and (\ref{e8}), we can compute $\Delta f$,
\allowdisplaybreaks
\begin{eqnarray*}
\Delta\frac{|H|^2}{\cos^{\frac{1}{b}}(\sqrt{ab}\theta)}
&=&\frac{1}{\cos^{\frac{1}{b}}(\sqrt{ab}\theta)}\Delta
|H|^2+|H|^2\Delta\frac{1}{\cos^{\frac{1}{b}}(\sqrt{ab}\theta)}+2\nabla
|H|^2\cdot\nabla\frac{1}{\cos^{\frac{1}{b}}(\sqrt{ab}\theta)}\\
&\geq&\cos^{-\frac{1}{b}}(\sqrt{ab}\theta) ((4-4\tilde b)|\nabla
|H||^2-\tilde a|H|^4-v^i\nabla_i |H|^2)\\
&&-|H|^2 v^i\nabla_i\cos^{-\frac{1}{b}}(\sqrt{ab}\theta)
+(\frac{a}{b}+a)\cos^{-\frac{1}{b}}(\sqrt{ab}\theta)\tan^2(\sqrt{ab}\theta)|H|^4
\\&& +a\cos^{-\frac{1}{b}}(\sqrt{ab}\theta)|H|^4
-2\frac{a}{b}\cos^{-\frac{1}{b}}(\sqrt{ab}\theta)\tan^2(\sqrt{ab}\theta)|H|^4
\\&&+2\cos^{\frac{1}{b}}(\sqrt{ab}\theta)\nabla f\cdot \nabla\frac{1}{\cos^{\frac{1}{b}}(\sqrt{ab}\theta)}
\\&\geq& \cos^{-\frac{1}{b}}(\sqrt{ab}\theta)\left[(4-4\tilde b)|\nabla
|H||^2-\tilde a|H|^4+(\frac{a}{b}+a)\tan^2{\sqrt{ab}\theta}|H|^4\right.\\
&&\left.+a|H|^4-2\frac{a}{b}\tan^2{\sqrt{ab}\theta}|H|^4\right]\\
&&-v^i\nabla_i f +2\cos^{\frac{1}{b}}(\sqrt{ab}\theta)\nabla
f\cdot \nabla\frac{1}{\cos^{\frac{1}{b}}(\sqrt{ab}\theta)}.
\end{eqnarray*}

Since $f$ is bounded, and the second fundamental form of $\Sigma$
is bounded which implies that the curvature of $\Sigma$ is
bounded, we can apply the generalized maximal principle (\cite{CY}
Theorem 3) to conclude that there is a sequence $\{x_k\}$ in $
\Sigma$, such that,
$$
\lim_{k\to\infty}f(x_k)=\sup_{x\in\Sigma}f(x),
$$
$$
\lim_{k\to\infty}|\nabla f(x_k)|=0,
$$
and
$$
\lim_{k\to\infty}\triangle f(x_k)\leq 0.
$$

It implies that, at $x_k$,
$$ \nabla
|H|^2=|H|^2\sqrt{\frac{a}{b}}\tan(\sqrt{ab}\theta)\nabla\theta
+o_k(1),
$$ where $o_k(1)\to 0$ as $k\to\infty$, thus we get that,
$$ |\nabla
|H||^2=\frac{a}{4b}\tan^2(\sqrt{ab}\theta)|H|^4+o_k(1).
$$
Putting this identity into the above inequality, we obtain that,
at point $x_k$,
\begin{eqnarray*}
\Delta\frac{|H|^2}{\cos^{\frac{1}{b}}(\sqrt{ab}\theta)} &\geq&
\cos^{-\frac{1}{b}}(\sqrt{ab}\theta)\left[\frac{a}{b}(1-\tilde b)
\tan^2(\sqrt{ab}\theta)|H|^4+(a-\tilde a)|H|^4\right.\\
&&\left.+(a-\frac{a}{b})\tan^2(\sqrt{ab}\theta)|H|^4\right] +o_k(1)\\
&\geq& \cos^{-\frac{1}{b}}(\sqrt{ab}\theta) \left[(a-\frac{\tilde
b}{b}a)\tan^2(\sqrt{ab}\theta)|H|^4 +(a-\tilde
a)|H|^4\right]+o_k(1).
\end{eqnarray*}
Set $b=\tilde b$ and $a>\tilde a$. Then at point $x_k$,
$$(a-\tilde a)\cos^{-\frac{1}{b}}(\sqrt{ab}\theta)|H|^4\leq o_k(1),
$$ which implies that $\lim_{k\to\infty}|H|^2(x_k)=0$. This is equivalent to
$\lim_{k\to\infty}f(x_k)=0$. That is not possible. Thus there is
no translating soliton with $\cos(\sqrt{ab}\theta)\geq\delta_0>0$.
This implies that $\inf_{\Sigma}\cos(\sqrt{ab}\theta)\leq 0$, i.e,
$\sup_{\Sigma}\sqrt{ab}|\theta|\geq\frac{\pi}{2}$. Now we estimate
$ab$. Since $\frac{2}{\varepsilon}=2+\delta$, so
$\delta=\frac{2-2\varepsilon}{\varepsilon}$
\begin{eqnarray*}
ab>\tilde a\tilde
b&=&\frac{2}{\varepsilon}(\frac{1}{2}+(1+\frac{1}{\delta})\frac{\varepsilon}{|T|-\varepsilon})
\\&=& \frac{2}{\varepsilon}(\frac{1}{2}+\frac{2-\varepsilon}{2-2\varepsilon}\frac{\varepsilon}{|T|-\varepsilon}) \\ &=&
\frac{1}{\varepsilon}+\frac{2-\varepsilon}{(1-\varepsilon)(|T|-\varepsilon)}.
\end{eqnarray*} Choose $\varepsilon=\frac{|T|}{|T|+1}$, then
$ab>2\frac{(|T|+1)^2}{|T|^2}.$ Therefore,
$\sup_{\Sigma}|\theta|>\frac{\sqrt{2}\pi}{4}\frac{|T|}{|T|+1}$.
This completes the proof.

\hfill Q. E. D.

Before proving the Main Theorem $1$, we need the following lemma.
\begin{lemma}\label{l1} Suppose that $\Sigma$ is a surface in
${\mathbf{R}}^4$, we have,
$$|\nabla\alpha|^2\leq |\overline\nabla
J_\Sigma|^2,$$ at the points where $\alpha$ is smooth.
\end{lemma}
{\it Proof.} In fact we can choose the local orthonormal frame of
$\{e_1, e_2, v_1, v_2\}$ on ${\mathbf{R}}^4$ along $\Sigma$ so
that $\omega$ takes the following form (cf. \cite{CL1}, \cite{CT},
\cite{CW} ),
$$\omega=\cos\alpha u_1\wedge u_2+\cos\alpha u_3\wedge
u_4+\sin\alpha u_1\wedge u_3-\sin\alpha u_2\wedge u_4
$$ where $\{u_1, u_2, u_3, u_4\}$ is the dual frame of $\{e_1, e_2, v_1,
v_2\}$, and
\begin{eqnarray*}
J =\left ( \begin{array}{clcr} 0 &\cos\alpha &\sin\alpha &0\\
-\cos\alpha &0 &0 &-\sin\alpha\\ -\sin\alpha &0 &0 &\cos\alpha
\\ 0 &\sin\alpha &-\cos\alpha &0\end{array}\right).
\end{eqnarray*}
Then $$\cos\alpha=\omega(e_1, e_2).$$ For the sake of simplicity,
we can assume the covariant derivatives of the orthonormal frame
$\{e_1, e_2, v_1, v_2\}$ satisfy
$$\nabla_{e_i}e_j=0.$$   We see that,
\begin{eqnarray*}
\nabla_1\cos\alpha &=&\omega(\overline\nabla_1 e_1,
e_2)+\omega(e_1,
\overline\nabla _1 e_2)\\
&=&h^\alpha_{11}\omega( v_\alpha, e_2)+h^\alpha_{12}\omega( e_1,
v_\alpha) \\&=&(h^2_{11}+h^1_{12})\sin\alpha,
\end{eqnarray*}
and
\begin{eqnarray*}
\nabla_2\cos\alpha &=&\omega(\overline\nabla_2 e_1,
e_2)+\omega(e_1,
\overline\nabla _2 e_2)\\
&=&h^\alpha_{21}\omega( v_\alpha, e_2)+h^\alpha_{22}\omega( e_1,
v_\alpha) \\&=&(h^1_{22}+h^2_{12})\sin\alpha,
\end{eqnarray*} Thus
$$|\nabla\cos\alpha|^2=(|h^2_{11}+h^1_{12}|^2+|h^1_{22}+h^2_{12}|^2)\sin^2\alpha
$$ i.e,
\begin{eqnarray*} |\nabla\alpha|^2 &=&(|h^2_{11}+h^1_{12}|^2+|h^1_{22}+h^2_{12}|^2)
\\&\leq& |\overline\nabla
J_{\Sigma}|^2 = |h_{11}^2+h_{12}^1|^2+|h^2_{21}+h^1_{22}|^2
+|h_{12}^2-h_{11}^1|^2+|h^2_{22}-h^1_{21}|^2  .
\end{eqnarray*}
\begin{remark} It is not hard to see that $|\nabla\alpha|^2=\frac{1}{2}|\overline\nabla
J_\Sigma|^2$ if $|H|^2=0$.
\end{remark}

Now we begin to prove Main Theorem $1$.

\vspace{.1in}

\noindent {\bf Main Theorem 1 } {\it Suppose that $\Sigma\in ST$
is a symplectic standard translating soliton with K\"ahler angle
$\alpha$, then $\sup_{\Sigma} |\alpha|>
\frac{\pi}{4}\frac{|T|}{|T|+1}$, where $T$ is the direction in
which the surface translates.}

\vspace{.1in}

{\it Proof.} We also prove it by contradiction. Assume that there
is a translating soliton with
$\cos(\sqrt{ab}\alpha)\geq\delta_0>0$, where $a, b$ are constants
which are determined later. We consider the function
$$f=\frac{|H|^2}{\cos^{\frac{1}{b}}(\sqrt{ab}\alpha)}.$$ From
(\ref{e1}) we see that, at the point where $\alpha$ is smooth,
$$\sin\alpha\Delta\alpha=(|\overline\nabla
J_{\Sigma}|^2-|\nabla\alpha|^2) \cos\alpha+v^i\nabla_i\cos\alpha.
$$ Then at the point where $\sin\alpha\neq 0$ (Note that $\alpha$ is smooth at the point
where $\sin\alpha\neq 0$), we have,
$$ \Delta\alpha=\frac{\cos\alpha}{\sin\alpha}(|\overline\nabla
J_{\Sigma}|^2-|\nabla\alpha|^2)-v^i\nabla_i\alpha.
$$
Thus using (\ref{e8}) we obtain that, \allowdisplaybreaks
\begin{eqnarray*}
\Delta\frac{|H|^2}{\cos^{\frac{1}{b}}(\sqrt{ab}\alpha)}
&=&\frac{1}{\cos^{\frac{1}{b}}(\sqrt{ab}\alpha)}\Delta
|H|^2+|H|^2\Delta\frac{1}{\cos^{\frac{1}{b}}(\sqrt{ab}\alpha)}+2\nabla
|H|^2\cdot\nabla\frac{1}{\cos^{\frac{1}{b}}(\sqrt{ab}\alpha)}\\
&\geq&\cos^{-\frac{1}{b}}(\sqrt{ab}\alpha) ((4-4\tilde b)|\nabla
|H||^2-\tilde a|H|^4-v^i\nabla_i |H|^2)\\
&&-|H|^2 v^i\nabla_i\cos^{-\frac{1}{b}}(\sqrt{ab}\alpha)
\\&&+\sqrt{\frac{a}{b}}\cos^{-\frac{1}{b}}(\sqrt{ab}\alpha)
\cot\alpha\tan(\sqrt{ab}\alpha)(|\overline\nabla
J_{\Sigma}|^2-|\nabla\alpha|^2)|H|^2
\\&&
+(\frac{a}{b}+a)\cos^{-\frac{1}{b}}(\sqrt{ab}\alpha)\tan^2(\sqrt{ab}\alpha)|H|^2
|\nabla\alpha|^2
\\&& +a\cos^{-\frac{1}{b}}(\sqrt{ab}\alpha)|H|^2|\nabla\alpha|^2
\\ &&-2\frac{a}{b}\cos^{-\frac{1}{b}}(\sqrt{ab}\alpha)\tan^2(\sqrt{ab}\alpha)|H|^2|\nabla\alpha|^2
\\&&+2\cos^{\frac{1}{b}}(\sqrt{ab}\alpha)\nabla f\cdot \nabla\frac{1}{\cos^{\frac{1}{b}}(\sqrt{ab}\alpha)}.
\end{eqnarray*}
Since $f$ is bounded, and the second fundamental form of $\Sigma$
is bounded which implies that the curvature of $\Sigma$ is
bounded, we can apply the generalized maximal principle (\cite{CY}
Theorem 3) to conclude that there is a sequence $\{x_k\}$ in $
\Sigma$, such that,
$$
\lim_{k\to\infty}f(x_k)=\sup_{x\in\Sigma}f(x),
$$
$$
\lim_{k\to\infty}|\nabla f(x_k)|=0,
$$
and
$$
\lim_{k\to\infty}\triangle f(x_k)\leq 0.
$$
 It implies that
$$ |\nabla
|H||^2=\frac{a}{4b}\tan^2(\sqrt{ab}\alpha)|\nabla\alpha|^2|H|^2+o_k(1),
$$ where $o_k(1)\to 0$ as $k\to\infty$.
By Lemma \ref{l1} and notice that
$\cot\alpha\tan(\sqrt{ab}\alpha)\leq \sqrt{ab}$, then we get that
at $x_k$, \allowdisplaybreaks
\begin{eqnarray*}
\Delta\frac{|H|^2}{\cos^{\frac{1}{b}}(\sqrt{ab}\alpha)} &\geq&
\cos^{-\frac{1}{b}}(\sqrt{ab}\alpha) \left[\frac{a}{b}(1-\tilde
b)\tan^2(\sqrt{ab}\theta)|\nabla\alpha|^2|H|^2\right.
\\ &&\left.+a(|\overline\nabla J_{\Sigma}|^2-|\nabla\alpha|^2)|H|^2\right.
\\&&\left. -\tilde a |H|^4+a|\nabla\alpha|^2 |H|^2\right. \\&&
\left. (a-\frac{a}{b})
\tan^2(\sqrt{ab}\theta)|\nabla\alpha|^2|H|^2\right] +o_k(1)\\
&\geq& \cos^{-\frac{1}{b}}(\sqrt{ab}\alpha) [(a-\frac{\tilde
b}{b}a)\tan^2(\sqrt{ab}\theta)|\nabla\alpha|^2|H|^2
\\ &&+(a|\overline\nabla J_{\Sigma}|^2-\tilde a |H|^2)|H|^2]+o_k(1).
\end{eqnarray*} Set $b=\tilde b$, $a=2\tilde  a$, by the same argument as Main Theorem 2, we get the
contradiction. Thus we choose $ab=2\tilde a\tilde
b>4\frac{(|T|+1)^2}{|T|^2}$. This completes the proof.

\hfill Q. E. D.

\section{A finite time blow-up example}

In this section, we construct a finite time blow-up example of
symplectic mean curvature flows in ${\bf C}^2$. Let $\gamma
(s)=(x(s),y(s))$ be a regular planar curve in ${\bf R}^2$ for
$s\in (-\infty,\infty)$. We consider the surface in ${\bf C}^2$
defined by
$$
F(s,\theta)=(x(s)e^{i\theta},y(s)e^{i\theta}):=\gamma(s)e^{i\theta},
$$
where $i=\sqrt{-1}$ and $\theta\in {\bf R}$.

Then
$$
F_s(s,\theta)=(x'(s)\cos\theta,x'(s)\sin\theta,y'(s)\cos\theta,y'(s)\sin\theta),
$$
$$
F_\theta(s,\theta)=(-x(s)\sin\theta,x(s)\cos\theta,-y(s)\sin\theta,y(s)\cos\theta).
$$

The normal vectors are
$$
\nu_s(s,\theta)=(-y'(s)\cos\theta,-y'(s)\sin\theta,x'(s)\cos\theta,x'(s)\sin\theta),
$$
$$
\nu_\theta(s,\theta)=(y(s)\sin\theta,-y(s)\cos\theta,-x(s)\sin\theta,x(s)\cos\theta).
$$

The induced metric is
$$
g_{ss}(s,\theta)=(x'(s))^2+(y'(s))^2,~~g_{\theta\theta}(s,\theta)=x^2(s)+y^2(s),~~g_{s\theta}=0.
$$

One checks that the K\"ahler angle is
$$
\cos\alpha
=\frac{(x(s)x'(s)+y(s)y'(s))}{\sqrt{(x^2(s)+y^2(s))((x'(s))^2+(y'(s))^2)}}.
$$

The second fundamental form is
$$
h_{ss}^s(s,\theta)=x'(s)y''(s)-x''(s)y'(s),~h_{\theta\theta}^s=x(s)y'(s)-x'(s)y(s),~
h_{s\theta}^s=h_{s\theta}^\theta=h_{\theta\theta}^\theta=0,
$$
and the mean curvature vector is
$$
H=\frac{1}{\sqrt{(x'(s))^2+(y'(s))^2}}(\frac{x'(s)y''(s)-y'(s)x''(s)}{(x'(s))^2+(y'(s))^2}
+\frac{x(s)y'(s)-x'(s)y(s)}{x^2(s)+y^2(s)})n(s)e^{i\theta},
$$
where $n(s)$ is the unit normal vector of $\gamma(s)$, and
$$
n(s)e^{i\theta}=\frac{1}{\sqrt{(x'(s))^2+(y'(s))^2}}\nu_s.
$$
These computations clearly yield the following proposition.
\begin{proposition} Let $r(s)=\sqrt{x^2(s)+y^2(s)}$. If the regular planar
curve $\gamma(s)=(x(s),y(s))$ satisfies that $r'(s)>0$, then the
surface in ${\bf C}^2$ defined by the map
$F(s,\theta)=(x(s)e^{i\theta},y(s)e^{i\theta})$ for $s,\theta\in
{\bf R}$ is symplectic. If $\gamma(t,s)=(x(t,s),y(t,s))$ satisfies
the curvature flow equation
\begin{equation}\label{curvef}
\frac{d}{dt}\gamma(t,s)=\frac{1}{\sqrt{(x')^2+(y')^2}}(\frac{x'y''-y'x''}{(x')^2+(y')^2}
+\frac{xy'-x'y}{x^2+y^2})n(t,s),
\end{equation}
where $()'=d/ds$, $()''=d^2/ds^2$, and $n(t,s)$ is the unit normal
vector of $\gamma(t,s)$, then
$F(s,\theta)=(x(t,s)e^{i\theta},y(t,s)e^{i\theta})$ satisfies the
mean curvature flow equation.
\end{proposition}

The curve flow (\ref{curvef}) first arises in the construction of
Lagrangian self-similar solutions to the mean curvature flow in
${\bf C}^n$ (c.f. \cite{A1}, \cite{A2}, \cite{GSSZ}, \cite{N1},
\cite{N2}).

The equation (\ref{curvef}) can also be written as
\begin{equation}\label{curvef1}
\frac{d}{dt}\gamma(t,s)=k(s,t)-\frac{\gamma^{\perp}(s,t)}{|\gamma
(s,t)|^2},
\end{equation}
where $k$ is the curvature of $\gamma$ and $\gamma^{\perp}$ is the
projection of the vector $\gamma$ on the orthogonal complement of
the tangent direction of $\gamma$.

\begin{theorem}
Let $\gamma_0(s)=\log(1+s)e^{is}$, $s\in [0,\infty)$. then the
curve flow (\ref{curvef}) with initial data $\gamma_0$ blows up at
a finite time.
\end{theorem}
{\it Proof:} Let $r_0(s)=\log(1+s)$. Then
$\gamma_0(s)=r_0(s)e^{is}$. It is proved in \cite{N1} (Lemma 4.3)
that the curve flow takes the form
$$
\gamma(t,s)=r(t,s)e^{is}.
$$
Computing directly, one sees that the equation (\ref{curvef}) is
deduced to the following form:
$$
\frac{d}{dt}\gamma(t,s)=\frac{3(r')^2-rr''+2r^2}{(r^2+(r')^2)^2}(-y',x'),
$$
because
$$
\frac{d}{dt}\gamma(t,s)=((\frac{d}{dt}r(t,s))\cos s,
(\frac{d}{dt}r(t,s))\sin s),
$$
the radius function $r(t,s)$ satisfies the equation
$$
\frac{d}{dt}r(t,s)=\frac{rr''-3(r')^2-2r^2}{r^3+r(r')^2}.
$$

We consider the weighted area $A(t)$ of the triangular region
$$
\{ue^{is}~|~0\leq s<\infty,~0\leq u\leq r(t,s)\},
$$
$$
A(t)=\frac{1}{2}\int_0^\infty\frac{r^2(t,s)}{(1+s)^2}ds.
$$
We have
\begin{eqnarray*}
\frac{d}{dt}A(t)&=&\int_0^\infty r(t,s)\frac{d}{dt}r(t,s)\frac{1}{(1+s)^2}ds\\
&=& \int_0^\infty\frac{1}{(1+s)^2}\frac{rdr'}{r^2+(r')^2}
-\int_0^\infty\frac{1}{(1+s)^2}\frac{(r')^2}{r^2+(r')^2}ds-2\\
&=& \int_0^\infty\frac{1}{(1+s)^2}\frac{d(r'/r)}{1+(r'/r)^2}-2\\
&\leq & 2\int_0^\infty\frac{1}{(1+s)^3}\arctan(r'/r)ds -2\\
&\leq & \frac{\pi}{2}-2<0.
\end{eqnarray*}

Therefore the curve flow must blow up at a finite time.

\hfill Q. E. D.

\end{document}